\documentclass[bib]{ba}

\usepackage{epsfig}
\usepackage{bayes}

\usepackage{amsmath,amsfonts,amsthm,amssymb,latexsym}
\usepackage{graphicx,picins}
\usepackage{natbib}

\newtheorem{example}{Example}[section]
\newcounter{exA}

\newcommand\dropcap\noindent

\begin{document}

\inserttype{article}
\author{Robert, C.P.}{%
  {\sc Christian P.~Robert}\\Universit\'e Paris-Dauphine, CEREMADE, and CREST, Paris
}	
\title[Keynes' Treatise on Probability]{{\it{\bfseries Reading Keynes'\\ {\em Treatise on Probability}}}}

\maketitle

\begin{abstract}
{\em A Treatise on Probability} was published by John Maynard Keynes in 1921. The {\em Treatise} contains a critical
assessment of the philosophical foundations of probability and of the statistical methodology at the time.  We review
the aspects of the book that are most related with statistics, avoiding uninteresting neophyte's forrays into
philosophical issues. In particular, we examine the arguments provided by Keynes again the Bayesian approach, as well as
the sketchy alternative of a return to Lexis' theory of analogies he proposes. Our conclusion is that the {\em Treatise}
is a scholarly piece of work looking at past advances rather than producing directions for the future.
%
\end{abstract}



\section{Introduction}
\newcommand\ToP{{\em{A Treatise on Probability}}}
%
\parpic[r]{\includegraphics[width=4.5cm]{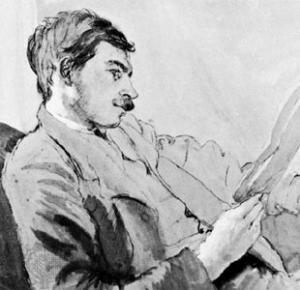}}

{\em{A Treatise on Probability}} is John Maynard Keynes' polishing of his 1907 and 1908 Cambridge Fellowship
dissertation ({\em The Principles of Probability}, submitted to King's College) into a book after an interruption due to
the war (for censorship reasons, as Keynes was then an advisor to the government). Although the author revised this
dissertation in 1914 (as mentioned by \citealp{aldrich:2008k}) and then in 1920 towards a general audience, the original
potential readers of \ToP~were therefore mostly local academics, among whom his Cambridge colleagues Edgeworth and Yule.
Despite Keynes' lasting interest in statistics, this is also his most significant publication in this field, since his
research focus had moved to economics by then \citep[as shown by, e.g.,][]{keynes:1919}. In contrast with the immense
influence Keynes exerted and still exerts in this latter field, and in agreement with the fact that the original version
was an internal dissertation, the impact of \ToP~was very limited, for reasons further discussed hereafter.

In this review, we consider the most relevant parts of the book, solely from a statistical perspective, avoiding the
outdated philosophical debates about the nature of probability and of induction that constitute most of the {\em
Treatise} and do not overlap with our personal interests. One must recall that this philosophical part on the
foundations of probability is a common feature in books of the period as shown by the first fifty pages of
\cite{jeffreys:1937} dedicated to ``direct probabilities".  Interestingly, Keynes favours a more subjective view of
probability as a degree of belief, \`a la de Finetti (\citeyear{definetti:1974}), while Jeffreys settles for a
mathematical definition that implies there is only one ``type" of probability. It is also worth reminding the reader at
this stage that Andrej Kolmogorov's book laying the axioms of modern probability was only published in
\citeyear{kolmogorov:1933}, since this explains why the concept of probability was still under debate at the time in
both mathematical and philosophical circles. 

As in the parallel review of \cite{jeffreys:1937} we undertook in \cite{robert:chopin:rousseau:2009}, there is no
attempt at drawing an history of statistics in this review, which is rather to be taken as a reflection of a modern
reader upon a piece of work written one century ago. For earlier historical aspects on the evolution of inverse
probability as a central piece of statistical thinking in the 19th Century, we refer the reader to the comprehensive
coverage in \cite{dale:1999} who, ironically, stops its range at Karl Pearson, i.e.~just before Keynes briefly entered
the statistical scene.  For a broader historical perspective on the development of statistics and the state of
statistics at the time of the {\em Treatise}, \cite{stigler:1986} undoubtedly remains the essential reference. %
\medskip

Before engaging upon this review, we point out that the {\em Treatise} has been previously assessed by \cite{stigler:2002},
who was similarly critical on the depth of the book, and by \cite{aldrich:2008k}, who produced an extensive and
scholarly survey on the impact (or lack thereof) of Keynes on the philosophical and statistical communities at the time.
The later includes in particular a detailed study of the reviews written on \ToP~by philosophers and statisticians of
the early 20th Century, including Ronald Fisher (\citeyear{fisher:1923}) and Harold \cite{jeffreys:1931}. 

\section{Contents of the {\em Treatise}}
    \begin{quote}{\em 
    `` A definition of {\em probability} is not possible, unless it contents us
    to define degrees of the probability-relation by reference to degrees of
    rational belief." {\em{\bfseries A Treatise on Probability}}, page 8.
    }\end{quote}

As clearly stated in the above quote, the proclaimed and ambitious goal of \ToP~is to establish a logical basis for
probability and of drawing a new ``constructive" approach for statistical induction. The extremely strong views
contained of the book, as well as the highly critical reassessments of past and (then) current authors, like Laplace and
Pearson, are reflecting upon the youth of the author and his earlier dispute with Karl Pearson on correlation. The
extensive coverage of the (statistical if not probabilistic) literature of the time and the comprehensive---if not always
insightful---discussion of most theories in competition shows the extent of the scholarly expertise of Keynes in
statistics. 

The {\em Treatise} comprises of 23 chapters, regrouped in five parts: 
\begin{itemize}
\item[I.] Fundamental Ideas;
\item[II.]  Fundamental Theorems;
\item[III.]  Induction and Analogy;
\item[IV.]  Some Philosophical Applications of Probability;
\item[V.]  The Foundations of Statistical Inference. 
\end{itemize}
The first part sets the logical and philosophical grounds for establishing a theory of probability, also touching upon
the Principle of Indifference treated below in Section \ref{sec:difR}. The second part is about probability axioms seen
from a mathematical logic perspective, although the mathematical depth is somewhat limited. This part also contains a
chapter on Inverse Probability, as discussed in Section \ref{sec:InvProb}. The third part is mostly philosophical and
discusses Humean induction with few connections with statistical inference. Part IV is a short metaphysical digression
on the meaning of randomness and its impact on conduct, completely unrelated with statistical inference.  The
statistical entries in the {\em{Treatise}} are mostly found in Part V, which covers convergence theorems (the "Law of
Great Numbers" and the "Theorems of Bernoulli, Poisson and Tchebycheff"), Bayesian inference and a call for a return to
the "Continental" principles laid by Lexis. As explored below, the amount of methodological innovation found in the book 
is extremely limited, in line with Keynes' own acknowledgement that he is {\em ``unlikely to get much further"}. 

\section{A restricted perspective}

From a statistician's viewpoint, the innovative aspects of the {\em Treatise} are quite limited in that the statistical
discourse remains at a highly rethorical---as opposed to methodological---level, drafting in vague terms the direction
for prospective followers that never materialised. While the {\em Treatise} presents both an historical
\citep{dale:1999} and a philosophical interest, from the perspectives both of Keynes' academic career and of the
foundations of statistics, there is no statistical advance to be found in it.  For instance, the {\em Treatise} is
missing the (then) current developments on a comprehensive theory of statistical tests, started with Karl Pearson's
$\chi^2$ and William Gosset's $t$ tests, and about to culminate in \cite{fisher:1925}.  Given the contents of the
soon-to-come major advances represented by not only Fisher's (\citeyear{fisher:1925}) {\em Statistical Methods}, but
also Jeffreys' (1939) and de Finetti's (1937, reprinted as 1974) homonymous {\em Theory of Probability}, the {\em
Treatise} does not stand the comparison as it fails to provide even a thorough treatment of the theory of statistics of
the time, if not proposing advances in this domain. 
This lack of innovative material, along with the harsh tone of a critic who had contributed so little to the field, may
explain why Keynes' incursion in probability and statistics did not have a lasting impact, since even those most
sympathetic to the book \citep{jeffreys:1931,lindley:1968} saw no practical nor methodological aspect to draw from and
praised aspects external to their own field \citep{aldrich:2008k}.  \citeauthor{stigler:2002} (2002, pp.~161-162)
similarly questions the worth of {\em{A Treatise on Probability}} as a mathematical and statistical work, with almost
sole focus {\em ``the binomial world"} and he considers the book unable to {\em ``carry the weight of a serious social
scientific investigation"}.

    \begin{quote}{\em 
    ``Statistical techniques tell us how to `count the cases'." {\em{\bfseries A Treatise on Probability}}, page 392.
    }\end{quote}

Keynes spends the major portion of the {\em Treatise} decrying a large part of the (then) current statistical practice
(first and foremost, Bayesian\footnote{\label{foot:soldier} It is worth pointing out that the denomination of `Bayesian'
appeared much later \citep{fienberg:2006}, Keynes resorting to the (then) current denomination of Inverse Probability
or, in a more derogatory way, {\em ``the Laplacian theory of `unknown probabilities''} (page 372). Referring to the
authoritative arguments of \cite{aldrich:2008k}, we stress that the papers of Keynes on statistics were definitely
Bayesian with most of his analysis being based on an uniform prior. He later started to worry about the influence of the
prior, leading to the {\em{Treatise}} and to its harsh criticism of inverse probability, even though some Bayesian
arguments remain in use within the book (see footnote \ref{noML}). This is quite in agreement with the practice of the
day, with statisticians mixing frequency and inverse probability arguments, Pearson and Fisher included, even though
earlier books like \cite{bertrand:1889} had pointed out the distinction between objective and subjective probabilities.}
statistics) as well as a majority of the past and (then) current statisticians, and in reproducing the arguments of
other (and more Continental) researchers, like Boole, Lexis, or von Kries.  (Again, this is in line with our argument
that the book is a scholarly and critical memoir rather than a innovative manifesto, even though the author aimed at a
broader impact on the statistics community at the time.) Furthermore, most of Part V deals with observation frequencies
(hence the quote at the top of this section) and their stabilisation. 

\parpic[r]{\includegraphics[width=4.5cm]{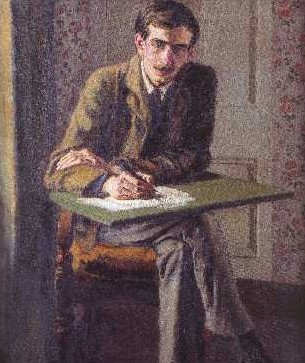}}
Quite curiously, when compared with, say, the much more modern treatise by
\cite{jeffreys:1931}\footnote{\cite{jeffreys:1922} reviewed the {\em Treatise} for {\em Nature}. His review was quite
benevolent, despite most of Keynes' perspectives on statistics being foreign to his own. This may explain why Part V
was mostly bypassed by Jeffreys' review. His comments in {\em Scientific Inference} (1931, pp.~222-224) about Keynes'
refusal to admit that probabilities were numbers, hence are comparable, and in {\em Theory of Probability} (1939, p.~25)
about Keynes' unwillingness to generalise axioms, more truthfully reflects Jeffreys' global opinion about the book.},
this book does not contain analyses of realistic datasets, except when criticising von Bortkiewicz's theory through the
Prussian cavalry horsekick data, which is customarily used for introducing Poisson modelling and is available in R as
the {\sf prussian} dataset \citep{cran}. This is somehow surprising when considering the main research field of Keynes,
namely Economics, where examples of considerable interest abound.  Instead, a very small number of (academic) examples
like the proportion of boys in births is recurrently discussed throughout the book. 

To be complete about the statistics contents of {\em{A Treatise on Probability}}, we note that Part II on {\em
Fundamental Theorems} also contains a chapter on the properties of various estimators of the mean in connection with the
distribution of the observations, although Keynes dismisses its importance by stating {\em ``It is without philosophical
interest and should probably be omitted by most readers"} (page 186). This chapter actually reproduces Keynes' only
genuine statistical paper, published in the {\em Journal of the Royal Statistical Society} in \citeyear{keynes:1911} on
the theory of averages (to be discussed in Section \ref{sec:1911}).
\medskip

\section{The low role of models}\label{sexymodels}

    \begin{quote}{\em 
    ``The knowledge of statistical theory, which is required for this, travels, I find, quite outside my knowledge."
      {\em{\bfseries Letter to Pearson, 1915}}.
    }\end{quote}

Throughout the book, Keynes holds both probability as a mathematical theory and probabilistic models as the basis for
statistics in very low regard, considering that unknown probability do not exist and that the reproducibility of
experiments is almost always questionable (as shown by the quote {\em `Some statistical frequencies are, with narrower
or wider limits, stable. But stable frequencies are not very common, and cannot be assumed lightly"}, page 336), apart
from urn models.  When adopting this type of reasoning, Keynes thus falls into what we call an ``ultra-conditioning
fallacy", namely that, the more covariates one conditions upon, the more different the individuals behave, a point of
view that goes against supporting statistical practice because there can be no frequency
stabilisation.\footnote{Discussing a similar point about Keynes, Stigler (1999, page 48) concludes that. with {\em
``this standard, it is difficult to conceive of any policy issue that can be investigated statistically"}.} For
instance, Keynes states that {\em ``where general statistics are available, the numerical probability which might be
derived from them is inapplicable because of the presence of additional knowledge with regards to the particular case"}
(page 29).  (He then goes on deriding Gibbon for his use of mortality tables when he should have called for a doctor!)
This shows the gap between the perspectives of Keynes and those of \cite{jeffreys:1931,jeffreys:1939} and de Finetti
(1937, reprinted as \citeyear{definetti:1974}), the later focussing on the exchangeability of events to derive the
existence of a common if unknown probability distribution. It is also the more surprising given the earlier works of
Pearson and Edgeworth in the 1880's that developed mathematical statistics towards a general theory of inference (see
\citealp{stigler:1986}, Chapters 9 and 10).  The above quote, given in \cite{stigler:1999} as a request from Keynes to
Pearson to help as examiner at the University of London, may however explain the reluctance of Keynes to engage in
deeper mathematics.

The use of particular sampling distributions (called {\em laws of errors}) in the reproduction of his 1911 paper on
averages (see Section \ref{sec:1911})
is not discussed in a modelling perspective but simply to back up the standard types of averages as maximum likelihood
estimators.\footnote{\label{noML}\cite{aldrich:2008k} argues quite convincingly that these are maximum a posteriori (MAP)
estimates corresponding to flat priors, making the term {\em most probable} more coherent. When stating the problem on
page 194, Keynes indeed invokes the {\em ``theorem of inverse probability"}, namely Bayes' theorem on the parameter {\em
finite} set. This theorem is also described in the 1911 paper. Again, this shows that part of Keynes' reasoning is still
grounded within the principles of inverse probability.} {\em ``The general evidence which justifies our assumption of
the particular law of errors which we do assume"} (page 195) is never discussed further in the {\em Treatise}. Assessing
the worth of a probability model against a dataset was not an issue for Keynes, although Pearson had earlier addressed
the problem.
\medskip

\section{Criticisms of frequentism}
    \begin{quote}{\em 
    ``The frequency theory, therefore, entirely fails to explain or justify the most important source of the most usual
      arguments in the field of probable inference." {\em{\bfseries A Treatise on Probability}}, page 108.
    }\end{quote}

Given Keynes' reluctance, mentioned above, to accept numerical probabilities, models and reproducibility, it is no
surprise that only extensive frequency stability is acceptable for him: {\em "The `Law of Large Numbers' is not at all a
good name for the principle that underlies Statistical Induction.  The `Stability of Statistical Frequencies' would be a
much better name for it."} (page 336). He has a strong {\em a priori} against the almost sure stabilisation of an iid
random sequence, especially when considering real data. This is historically intriguing given the derivations by
Bernoulli, de Moivre, and Laplace of the Law of Large Numbers more than a century earlier.

    \begin{quote}{\em 
    ``Some statistical frequencies are, with narrower or wider limits, stable. But stable frequencies are not
      very common, and cannot be assumed lightly." {\em{\bfseries A Treatise on Probability}}, page 336.
    }\end{quote}

\parpic[r]{\includegraphics[width=4.5cm]{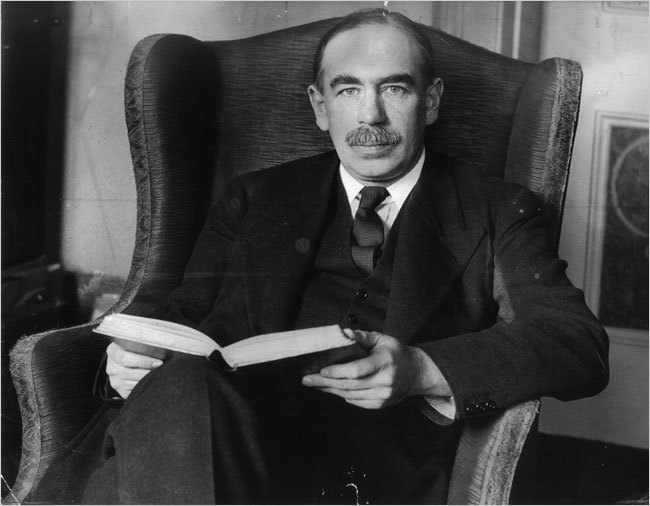}}
The criticism of the Central Limit Theorem (CLT), called Bernoulli's Theorem in the book,\footnote{Bernoulli's Theorem
is historically the weak Law of Large Numbers but Keynes presents this result in conjunction with (a) a description of
the binomial $\mathcal{B}(n,p)$ distribution and (b) the normal (CLT) approximation to the binomial cdf, a result he
calls Stirling's theorem. While Edgeworth had a clear influence on Keynes in Cambridge, his expansions providing a
better approximation than the CLT are not mentioned \citep{hall:1992}.} that is found in Chapter XXIX is rather curious, in that it
confuses model probabilities $p$ with probability estimates $p^\prime$ (the identical notation being an indicator of
this confusion).\footnote{Jeffreys (1931, p.~224) stresses that {\em ``Keynes' postulate might fit the assigned
probabilities instead of the true probabilities".}} For instance, on page 343, Keynes criticises the use of the Central
Limit Theorem (CLT) for the Bernoulli distribution $\mathcal{B}(1/2)$ and a coin tossing experiment as, when {\em
``heads fall at every one of the first 999 tosses, it becomes reasonable to estimate the probability of heads at much
more than $1/2$"}.  This argument is therefore confusing the probability model $\mathcal{B}(p)$ with the estimation
problem. Keynes' inability to recognise the distinction may stem from his reluctance to use unknown probabilities such
as $p$.  Similarly, on pages 349-350, when considering the proportion of male births, an example dating back to Laplace,
Keynes states that the probability of having $n$ males births in a row is not $p^n$, if $p$ is the probability of a single 
male birth, but 
$$
\dfrac{r}{s}\cdot\dfrac{r+1}{s+1}\cdot\dfrac{r+2}{s+2}\cdots\dfrac{r+n-1}{s+n-1}
$$
if $s$ is the number of births observed so far, and $r$ the number of male births. The later is a sequential
construction based on individual estimates for each new observation, neither a true (predictive) probability nor a
genuine plug-in estimate. Most interestingly, under a flat prior on $p$, the predictive (marginal) probability of seeing $n$ male
births in a row is
\begin{align*}
\int_0^1 p^n\,\dfrac{(s+1)!}{r!(s-r)!}\,p^r(1-p)^{s-r}\,\text{d}p 
&= \dfrac{(s+1)!}{r!(s-r)!}\,\dfrac{(r+n)!(s-r)!}{(n+s+1)!}\\
&= \dfrac{(r+1)\cdots(r+n)}{(s+2)\ldots(s+n+1)}\,r.
\end{align*}
We thus come to the conclusion that Keynes' solution corresponds to using Haldane's (\citeyear{haldane:1932}) prior,
$\pi(p)=1/(p(1-p)$, whose impropriety difficulties \citep{robert:2004} were not an issue at the time and even later in
\cite{jeffreys:1939}.

    \begin{quote}{\em 
    ``It seldom happens that we can apply Bernoulli's theorem to a long series of
      natural events." {\em{\bfseries A Treatise on Probability}}, page 343.
    }\end{quote}

That Keynes concludes that Bernoulli's Theorem (a simple version of the CLT in the binomial case) does not hold exactly
in this setting is clearly inappropriate. When considering that {\em ``knowledge of the result of one trial is capable
of influencing the probability at the next"}, he is confusing the ``true" probability with the estimated one.  The same
criticism applies to Keynes' remark that {\em ``a knowledge of some members of a population may give us a clue to the
general character of the population in question."} (page 346), a remark that bears witness to Keynes' skepticism about
the relevance of probabilistic models. From a Bayesian perspective, it appears that Keynes mixes sampling distributions
with marginal distributions, as in the latent variable example of page 346 dealing with observations from
$\mathcal{B}(p)$ when $p\in\{p_1,\ldots,p_k\}$: the observations become dependent when integrating out $p$. The
statement {\em ``if we knew the real value of the quantity, the different measurements of it would be independent"}
(page 195) may be understood under this light, even though it is a risky extrapolation given both the book stance on
Bayesian statistics and the lack of evidence Keynes mastered this type of mathematical techniques (shown by the quote at
the entry of Section \ref{sexymodels}).  
\medskip

\section{Keynes' views on  statistical inference}\label{sub:Kiew}

In the continuation of the quote from page 392 given above, the {\em Treatise} argues most vigorously against 
mathematical statistics by stating that the purpose of statistics ought to be {\em strictly limited to preparing
the numerical aspects of our material in an intelligible form}. Keynes thus separates inference ({\em the usual
inductive methods}) from statistics and clearly shows his skepticism about extending statistics beyond a descriptive
tool.

    \begin{quote}{\em 
    ``The statistician is less concerned to discover the precise conditions in which a description can be legitimately extended
    by induction." {\em{\bfseries A Treatise on Probability}}, page 327.
    }\end{quote}

The focus of statistical inference as described in Part V is reduced to a probability assessment: {\em ``In the first type
of argument we seek to infer an unknown statistical frequency from an {\em \`a priori} probability. In the second type
we are engaged on the inverse operation, and seek to base the calculation of a probability on an observed statistical
frequency. In the second type we seek to pass from an observed statistical frequency, not merely to the probability of
an individual occurrence, but to the probable value of other unknown statistical frequencies"} (page 331). This is
actually rather surprising given the overall negative tone of {\em{A Treatise}} about probability theory.

\parpic[r]{\includegraphics[width=4.5cm]{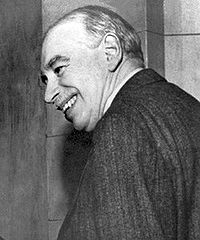}}
The first item above is a probabilistic issue and is treated as such in Chapter XXIX, which covers both the normal and
Poisson limit theorems, as well as \v Ceby\v sev's inequality.  Further criticisms of ``Bernoulli's Theorem" found in
this chapter are limited to the fact that finding independent and identically distributed (i.i.d.) replications is a
condition that is {\em ``seldom fulfilled"} (page 342). The 1901 proof by Liapounov of the CLT for general independent
random variables is not mentioned in Keynes' book and was presumably unknown to the author. Instead, he refers to
Poisson for a series of independent random variables with different distributions, warning that {\em ``it is important
not to exaggerate the degree to which Poisson's method has extended the application of Bernoulli's results"}  (page
346). Although \v Ceby\v sev's inequality has had very little impact on statistical practice, except when constructing
conservative confidence intervals, Keynes is clearly impressed by the result (of which he provides a very convoluted
proof on pages 353-355) and he concludes---rather unfairly since Laplace wrote one century before---that {\em 
``Laplacian mathematics is really obsolete and should be replaced by
the very beautiful work which we owe to these Russians"} (page 355). Chapter XXIX terminates with an interesting section
on simulation experiments aiming at an empirical verification of the CLT, although Keynes' conclusion on a very long
dice experiment is that, given that the frequencies do not match up {\em ``what theory would predict"} (page 363), the
dice used in this experiment was quite irregular (or maybe worn out by the 20,000 tosses!)

    \begin{quote}{\em 
    ``I do not believe that there is any direct and simple method by which we can make the transition from an
    observed numerical frequency to a numerical measure of probability." {\em{\bfseries A Treatise on Probability}}, page 367.
    }\end{quote}

As illustrated by the above quote and discussed in the next section, Keynes does not consider Laplace's (i.e.~the
Bayesian) approach to be logically valid and he similarly criticises normal approximations \`a la Ber\-noul\-li, seeing both
as {\em ``mathematical charlatanery"} (page 367)! Even the (maximum likelihood) solution of estimating $p$ with the
frequency $x/n$ when $x\sim\mathcal{B}(n,p)$ does not satisfy him (as being {\em ``incapable of a proof"}, page 371).
Note that the maximum likelihood estimator is called the {\em most probable value} throughout the book, in concordance
with the current denomination at the time \citep{hald:1999}, without Keynes objecting to its Bayesian flavour.
Obviously, given that he wrote the main part of the book before the war, he could not have used Fisher's denomination of
maximum likelihood estimation since its introduction dates from \citeyear{fisher:1922}.\footnote{Both "Fisher" and
"estimation" entries are missing from the index of the {\em Treatise}.}

The method of least squares is also heavily attacked in Chapter XVII as {\em ``surrounded by an unnecessary air of
mystery"} (page 209), while conceding on the next page that it exactly corresponds to assuming the normal distribution
on observations (a fact that is not correct either).  Once again, Keynes is missing the recent developments of Pearson
and Keynes on the estimation of regression coefficients, following the publication in \citeyear{galton:1889} of {\em
Natural Inheritance} by Francis Galton and his discovery of regression ({\em ``one of the most attractive triumphs in the
history of statistics"}, according to Stigler, 1999, page 186.). Galton is only quoted twice in the {\em Treatise} and
for marginal reasons, while regression does not appear at all.

\section{On the principal averages}\label{sec:1911}

Chapter XVII reproduces Keynes' 1911 paper in the {\em Journal of the Royal Statistical Society} on the characterisation
of the distributions leading to specific standard averages as MAPs under a flat prior, i.e.~modern MLEs, which means
obtaining classes of densities for which the MLEs are the arithmetic, the geometric and the harmonic averages, and the
median, respectively.  The earlier decision-theoretic justifications of the arithmetic mean by Laplace and Gau\ss~are
derided as depending on {\em ``doubtful and arbitrary assumptions"} (page 206), while the lack of reparameterisation
invariance of the arithmetic average as MLE is clearly stated (on page 208). This classification of standard averages as
MLEs is more of a technical exercise than of true methodological relevance, because the classification of distributions
({\em ``laws"}) that give the arithmetic, geometric, harmonic mean or the median as MLEs is obviously
parameterisation-dependent, a fact later noted by Keynes but omitted at this stage despite his criticism of Laplace's
principle on the same ground. The derivation of the densities $f(x,\theta)$ of the distributions is based on the
condition that the likelihood equation
$$
\sum_{i=1}^n \dfrac{\partial}{\partial\theta}\log f(y_i,\theta) = 0
$$
is satisfied for one of the four empirical averages, using differential calculus despite the fact that Keynes earlier
derived (on page 194) Bayes' theorem by assuming the parameter space to be discrete.\footnote{Keynes notes that {\em
``differentiation assumes that the possible values of $y$ {\em [meaning $\theta$ in our notations]} are so numerous and
so uniformly distributed that we may regard them as continuous"} (page 196).} Under regularity assumptions, in the case
of the arithmetic mean, this leads to the family of distributions
$$
f(x,\theta) = \exp\left\{ \phi^\prime(\theta) (x-\theta) - \phi(\theta) + \psi(x) \right\}\,,
$$
where $\phi$ and $\psi$ are arbitrary functions such that $\phi$ is twice differentiable and $f(x,\theta)$ is a density in $x$, meaning that
$$
\phi(\theta) = \log \int \exp \left\{ \phi^\prime(\theta) (x-\theta) + \psi(x)\right\}\, \text{d}x\,,
$$
a constraint missed by Keynes. (The same argument is reproduced in Jeffreys, 1939, page 167.) 

While we cannot judge of the level of novelty in Keynes' derivation with respect to earlier works, this derivation
interestingly produces a generic form of unidimensional exponential family, twenty-five years before their rederivation by
\cite{darmois:1935}, \cite{pitman:1936} and \cite{koopman:1936} as characterising distributions with sufficient
statistics of constant dimensions. The derivation of distributions for which the geometric and the harmonic means
are MLEs then follows by a change of variables, $y=\log x,\,\lambda=\log \theta$ and $y=1/x,\,\lambda=1/\theta$,
respectively. In those different derivations, the normalisation issue is treated quite off-handedly by Keynes, witness
the function
$$
f(x,\theta) = A \left( \dfrac{\theta}{x} \right)^{k\theta} e^{-k\theta}
$$
at the bottom of page 198, which is not integrable {\em per se}. Similarly, the derivation of the log-normal density on page 199 is missing
the Jacobian factor $1/x$ (or $1/y_q$ in Keynes' notations) and the same problem arises for the inverse-normal density, which should be
$$
f(x,\theta) = A e^{-k^2(x-\theta)^2/\theta^2 x^2} \dfrac{1}{x^2}\,,
$$
instead of $A\exp k^2(\theta-x)^2/x$ (page 200). At last, the derivation of the distributions producing the median as
MLE is rather dubious because it does not seem to account for the non-differentiability of the absolute distance in 
every point of the sample. Furthermore, Keynes' general solution
$$
f(x,\theta) = A \exp \left\{ \displaystyle{\int \dfrac{y-\lambda}{|y-\lambda|}\,\phi^{\prime\prime}(\lambda)\,\text{d}\lambda +\psi(x) }\right\}\,,
$$
where the integral is interpreted as an anti-derivative, is such that the recovery of Laplace's distribution, 
$f(x,\theta)\propto \exp-k^2|x-\theta|$ involves setting (page 201)
$$
\psi(x) = \dfrac{\theta-x}{|\theta-x|}\,k^2x\,,
$$
hence making $\psi$ dependent on $\theta$ as well. In his summary (pages 204-205), Keynes (a) reintroduces a constant $A$ 
for the normalisation of the density in the case of the arithmetic mean and (b) produces
$$
f(x,\theta) = A \exp \left\{ \phi^\prime(\theta)\,\dfrac{\theta-x}{|\theta-x|}+\psi(x)\right\}
$$
in the case of the median. This later form is equally puzzling because the ratio in the exponential is equal to the sign of $x-\theta$, leading
to a possibly different weighting of $\exp \psi(x)$ when $x<\theta$ and when $x>\theta$.
\medskip

\section{A reactionary proposal}

After an extensive criticism of the methods of the time and of the use of mathematical models as a basis for
statistical inference (see Section \ref{sexymodels}), Keynes concludes the {\em Treatise} with a defence of the
method advocated by the late Lexis (who died in 1914), at the very moment \cite{fisher:1925} was defining statistics as 
``mathematics applied to data". As analysed by Aldrich (\citeyear{aldrich:2008k}, Section 5), the defence is paired with Keynes' attempt to
link Lexis' theory to his own principles of analogy in induction, as advanced in Part III of the {\em Treatise}. The
following quote indicates why the attempt failed.

\begin{quote}
{\em ``I have experienced exceptional difficulty, as the reader may discover for himself in the following pages, both in
clearing up my own mind about it and in expounding my conclusions precisely and intelligently."} {{\bfseries A Treatise
on Probability}}, page 409.  
\end{quote}

\parpic[l]{\includegraphics[width=4.5cm]{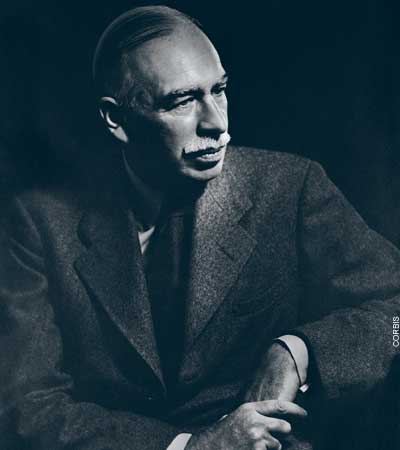}}
When considered from a modern perspective, Chapters XXXII and XXXIII advocate a very empirical approach to statistics
(which, in an anachronistic way, prefigures bootstrap), namely to derive the stability of a probability estimate by
subdividing a {\em series} into a large enough number of {\em subs Eries} in order to assess the variability of the
estimate or to spot heterogeneity. Keynes associates this approach with Lexis and appears quite supportive of the
latter, even though he comments that {\em ``Lexis has not pushed his analysis far enough"} (page 401), before
complaining about von Bortkiewicz, {\em ``preferring algebra to earth"} (page 404).  As highlighted by the above quote,
the {\em Treatise} faces difficulties in building a general theory around this approach and the description of the
mechanism for dividing the series remains unclear throughout the chapters, since it seems to depend on covariates. For
instance, the sentence {\em ``all conceivable resolutions into partial groups"} (page 395) is to be opposed to breaking
{\em ``statistical material into groups by date, place, and any other characteristic which our generalisation proposes
to treat as irrelevant"} (page 397). The model thus constructed has a mixture flavour when the groups are made per chance,
or a hierarchical one otherwise.  Indeed, the description of {\em ``the probability $p$ for the group made up as follows"}
(page 395) $$ p = \dfrac{z_1}{z}p_1+\dfrac{z_2}{z}p_2+\ldots $$ clearly corresponds to a mixture, the $z_i$'s being the
component sizes. 

In any case, the description of Lexis' theory sums up as testing for variations between groups, i.e. by exposing a
possible extra-binomial---called {\em supra-normal} by Keynes---variation. Keynes also mentions the possibility of an
insufficient variation---the {\em subnormal} case---is attributed to dependence in the data, which {\em ``cannot be
handled by purely statistical methods"} (page 399)\footnote{The whole book considers handling dependent observations an
impossible task, despite Markov's introduction of Markov chains a few years earlier. As pointed out by Aldrich (2008),
the lack of feedback from Yule in {\em{A Treatise on Probability}} is apparent from the pessimistic views of the author
about dependent series, despite the proximity of the authors in Cambridge, since Yule had already engaged into building
a statistical analysis of time series.}. A modern accounting of Lexis' procedure for testing stability and of why the
author ``failed" (and is now largely forgotten) is given in Stigler (1986, Chapter 6), who adds on page 238 that Keynes,
as one of the few followers of Lexis, missed the point that {\em ``simple urns models were insufficiently rich to
support the needs of a modern statistical analysis"}.

    \begin{quote}{\em 
    ``Statistical induction is not really about the particular instance at all but a {\em series}."
    {\em{\bfseries A Treatise on Probability}}, page 411.
    }\end{quote}

As the final chapter, Chapter XXXIII is Keynes' last attempt at defending his own views about a constructive theory of
statistical inference. However, it mostly sounds like rephrasing Lexis' views, Keynes' main point being that one should
work with {\em ``{\em series of series} of instances"} (page 407) in order to check for the stability of the assumed
model. The point made in the above quote is valid at face value but the attempt at checking that all subdivisions of a
dataset show the same variability ({\em ``until a {\em prima facie} case has been established for the existence of a
stable probable frequency, we have but a flimsy basis for any statistical induction"}, page 415) is doomed when pushed
to its extreme division of the data into individual observations. Furthermore, we again stress that the {\em Treatise}
never explicitly derives a testing methodology in the sense of Gosset or of Fisher\footnote{Ronald Fisher also reviewed
Keynes' book in \citeyear{fisher:1923}, concluding at the uselessness of Keynes' perspective on statistics as described
in \cite{aldrich:2008}.}, despite mentions made of {\em ``significant stability"} (pages 408 and 415).  When discussing
Lexis' dispersion in Chapter XXXII, Keynes refers to a case when {\em ``the dispersion conformed approximately to the
(...) normal law of error"} (page 358), but, again, no entry is found on the contemporary Student's $t$ tests or
Pearson's $\chi^2$ tests.  The earlier criticisms of Keynes' about the extension of an observed model to future
occurrences also apply in this setting, a fact acknowledged by the author: {\em ``it is not conclusive and I must leave
to others its more exact elucidation"} (page 419).  Besides, the assessment of stability is not detailed and, while it
seems to be based on normal approximations (to the binomial), the facts that the {\em same} data is used repeatedly and
thus that the test statistics are dependent appear to have been overlooked by Keynes.
\medskip

\section{Inverse Probability}\label{sec:InvProb}

As already discussed in footnotes \ref{foot:soldier} and \ref{noML}, the foundations of statistics were not sufficiently
settled at the time Keynes wrote his book to allow for a clear distinction between frequentist and Bayesian
philosophies. The choice of prior distributions had already come under attacks in the books of Chrystal, Venn and
Bertrand, but the alternative construction of a non-Bayesian setting would have to wait a few more
years for Fisher's (\citeyear{fisher:1925}) new perspective.

    \begin{quote}{\em 
    ``Bayes' enunciation is strictly correct and its method of arriving at it shows its true logical connection
    with more fundamental principles, whereas Laplace's enunciation gives it the appearance of a {\em new} principle
    specially introduced for the solution of causal problems." {\em{\bfseries A Treatise on Probability}}, page 175.
    }\end{quote}

When discussing the history of Bayes' theorem in Chapter XVI, Keynes considers---as shown by the above quote---that only
Bayes got his proof right and that subsequent writers, first and foremost Laplace, muddled the issue (except for
Markov)! While the author of the {\em Treatise} rightly separates the mathematical result represented by Bayes' theorem
from its use in statistical inference, Keynes misses the fact that Laplace independently derived Bayes' theorem from a
purely mathematical perspective, before applying (much later) inverse probability principles in statistical problems.
(Misunderstanding Bayes' theorem with not-yet Bayesian statistics seemed to be quite common at the time since, as
reported in \citealp{stigler:1999}, Karl Pearson equates Bayes' theorem with Laplace's Principle of Non-Sufficient
Reason covered below.)

An interesting discussion in the {\em Treatise} revolves around the (obvious) fact that the prior probabilities of the
different causes should be taken into account ({\em ``the necessity in general of taking into account the {\em \`a
priori} probability of the different causes"}, page 178). But one argument sheds light on the difficulty Keynes had with
the updating of probabilities, as mentioned in the paragraph about the CLT: {\em ``how do we know that the possibilities
admissible {\em \`a posteriori}  are still, as they were assumed to be {\em \`a priori}, equal possibilities} (page
176).\footnote{Note the accents used in {\em \`a priori}  and {\em \`a posteriori}, although there are no accent in
Latin. They may have stemmed from the way French writers first used those terms, even though {\em \`a posteriori} is
also found in Jakob Bernoulli's {\em Ars Conjectandi}...  The accent has vanished by the time of \cite{jeffreys:1931}.}
This section considers the specific arguments Keynes advanced against Bayesian principles. (We note again that the
statistical practice had the time had both frequentist and Bayesian, i.e.~sampling and posterior, arguments mostly mixed
in its arguments, as detailed in \citealp{aldrich:2008k}.)
\medskip

\subsection{Against the Principle of Indifference}\label{sec:difR}
    \begin{quote}{\em 
    ``My criticism will be purely destructive and I will not attempt to indicate my own way out of
    the difficulties." {\em{\bfseries A Treatise on Probability}}, page 42.
    }\end{quote}

The {\em Principle of Indifference} is Keynes' renaming of the {\em Principle of Non-Sufficient Reason} advocated by
Laplace and his followers for using (possibly improper) uniform prior distributions. Following the above preamble,
Keynes (rightly) shows the inconsistency of this approach under (a) a refinement of the available alternative (pages
42-43) and (b) a non-linear reparameterisation of the model (page 45), the example being the change from $\nu$ into
$1/\nu$. An extension of this argument on page 47 discusses the dependence of the uniform distribution on the dominating
measure (although the book obviously does not dally with a measure theory not yet finalised at the time and anyway
beyond Keynes' reach), as illustrated by Bertrand's paradox.\footnote{We note however that the general setup of modern
measure theory had already been given by Henri Lebesgue in 1903 in {\em Annali di Mathematica}.} This paradox points out
the lack of meaning of a ``random chord" of a circle without a proper probability structure and is reanalysed in Jaynes
(\citeyear{jaynes:2003}, page 386) from an objective Bayes perspective, where the author defends the maximum invariance
principle. In the following and less convincing paragraphs of Chapter IV, Keynes finds defaults with basic game examples
(including the Monty Hall problem) where again the equidistribution depends on the reference measure. 

\begin{quote}{\em 
``Who could suppose that the probability of a purely hypothetical event, of whatever complexity
(...), and which has failed to occur on the one occasion on which the hypothetical conditions were fulfilled is no less
than $1/3$?." {\em{\bfseries A Treatise on Probability}}, page 378.  
}\end{quote} 

Similar arguments are advanced in
Chapter XXX when debating about Laplace's law of succession.\footnote{Although no black swan glides in, the section
contains the obligatory example of the probability of the sun rising tomorrow that found in almost every treatise on
induction since Hume \citep{taleb:2004}.} Those are standard criticisms found for instance in the earlier
\cite{bertrand:1889}. Namely, putting a uniform distribution on all possible alternatives is not coherent given that a
subdivision of an alternative into further cases modifies the uniform prior. And, furthermore, a non-linear
reparameterisation of a probability $p$ into $q=p^n$ fails to carry uniformity from $p$ to $q$.  In concordance with the
spirit of the time \citep{lhoste:1923,broemeling:2003}, the debate about whether or not the Principle of Indifference
holds makes some sense, as shown by the subsequent defence by Jeffreys, but it does not hold much appeal nowadays
because priors are recognised as reference tools for handling data rather than expressions of truth or of ``objective
probabilities".

Chapter XXXI debates on the inversion of Bernoulli's Theorem, a notion that we interpret as Bayes formula applied to the Gaussian approximation
to the distribution of an empirical frequency: on page 387
$$
\dfrac{f(q^\prime)|h \cdot f(q)}{\sum f(q^\prime)|h \cdot f(q)}\,, 
$$
apparently meaning
$$
\dfrac{f(x|\theta) \pi(\theta)}{\int f(x|\theta) \pi(\theta)\,\text{d}\theta}
$$
in modern notations, is associated with the statement that {\em ``all the terms can be determined numerically by
Bernoulli's Theorem"}. Since this representation is somehow based on a flat reference prior (although the formula at the
bottom of page 386 which seems to involve two distributions on the parameter $\theta$ is incomprehensible), and thus on
the Principle of Indifference, it is rejected by Keynes who cannot see a {\em ``justification for the assumption that
all possible values of $q$ are {\em \`a priori} equally likely"} (page 387).  \medskip

\subsection{Against probabilising the unknown}
    \begin{quote}{\em 
    ``Laplace's theory requires the employment of both of two inconsistent methods." {\em{\bfseries A Treatise on Probability}}, page 372.
    }\end{quote}

The criticism of Bayesian (Laplacian's) techniques goes further than the rather standard debate about the choice of the
prior. For Keynes, adopting a perspective that unknown probabilities can be modelled as random variables is beyond
logical reasoning. Because an {\em unknown} probability is indeterminate, Keynes considers that {\em ``there is {\em no}
such value"} (page 373).\footnote{A more positive perspective (see, e.g., \citealp{brady:2004}) is to consider that
Keynes' stance prefigures the theory of imprecise probabilities \`a la Dempster--Schafer (see, e.g.,
\citealp{walley:1991}), as for instance when he states that {\em ``many probabilities can be placed {\em between}
numerical limits"} (page 160), but, to us, this mostly shows the same confusion between (interval) estimates and true
probabilities found elsewhere in the book. The notion of replacing (pointwise) probabilities by intervals in the short
Chapter XV is attributed to Boole---with another barb in the footnote of page 161---and it does not seem to be set to
any implementable version in the statistical inference section (Part V).} Therefore, the Bayesian notion of setting a
probability distribution over the unit interval is both illogical and impractical, since {\em ``if a probability is
unknown, surely the probability, relative to the same evidence, of this probability has a given value, is also unknown"}
(page 373). Keynes then argues that, if the hyperprior probability is unknown, it should also be endowed with its own
probability measure, inducing {\em ``an infinite regress"} (page 373).
\medskip

\section{Conclusion}

In conclusion, while Keynes' early interest in Probability and in Statistics is unarguable, {\em{A Treatise on
Probability}} could not have made a lasting contribution to Statistics, even from an historical perspective, given the
immense developments taking place in Statistics at the turn of the Century or in the neighbouring decades. The {\em
Treatise} appears in the end as a scholarly exercise focussing on past books and lacking a vision of developments that
would have made Keynes a statistician of his time, while the aggressive tone adopted towards most of the writers quoted
in the book is undeserved when comparing the achievements of both camps. It is therefore no surprise the book has had no
influence on the probability and statistics communities: it would make no sense to advise students in the field to put
aside major treatises to ponder through {\em{A Treatise on Probability}} as, to adopt Fisher's (1922) harsh but still
relevant words, {\em ``they would be turned away, some in disgust, and most in ignorance, from one of the most promising
branches of mathematics."}

\section*{Acknowledgements}
The author's research is partly supported by the Agence Nationale de la Recherche (ANR, 212, rue de Bercy 75012 Paris)
through the 2007--2010 grant ANR-07-BLAN-0237 ``SPBayes". The first draft of this paper was written during the
conference on Frontiers of Statistical Decision Making and Bayesian Analysis in San Antonio, Texas, March 17-20, held in
honour of Jim Berger's 60th birthday, and the author would like to dedicate this review to him in conjunction with this
event. He is also grateful to Eric S\'er\'e for his confirmation of the arithmetic mean distribution classification
found in Chapter XVII of the {\em Treatise}. Detailed and constructive comments from a referee greatly helped in
preparing the revision of the paper. This paper was composed using the {\sf ba.cls} macros from the International Society
for Bayesian Analysis.

\renewcommand{\bibsection}{\section*{References}}


\end{document}